
\input epsf.tex
\input amssym.def
\input amssym
\magnification=1100
\baselineskip = 0.23truein
\lineskiplimit = 0.01truein
\lineskip = 0.01truein
\vsize = 8.7truein
\voffset = 0.1truein
\parskip = 0.10truein
\parindent = 0.3truein
\settabs 12 \columns
\hsize = 5.8truein
\hoffset = 0.2truein

\setbox\strutbox=\hbox{%
\vrule height .708\baselineskip
depth .292\baselineskip
width 0pt}
\font\caps=cmcsc10
\font\bigtenrm=cmr10 at 14pt

\def\sqr#1#2{{\vcenter{\vbox{\hrule height.#2pt
\hbox{\vrule width.#2pt height#1pt \kern#1pt
\vrule width.#2pt}
\hrule height.#2pt}}}}
\def\square{\mathchoice\sqr46\sqr46\sqr{3.1}6\sqr{2.3}4}

\centerline{\bigtenrm SPECTRAL GEOMETRY, LINK COMPLEMENTS}
\centerline{\bigtenrm AND SURGERY DIAGRAMS}
\tenrm
\vskip 14pt
\centerline{MARC LACKENBY\footnote*{Supported by an EPSRC Advanced Research Fellowship}}
\vskip 18pt

\centerline{\caps 1. Introduction}
\vskip 6pt

The spectrum of the Laplacian on a Riemannian manifold $M$ has
been the focus of an enormous amount of study. Of particular
importance is $\lambda_1(M)$, which is the infimum of ${\rm Spectrum}(M) - \{ 0 \}$.
This contains a good deal of geometric information
about $M$; in particular, it is closely related to a geometric
quantity, the Cheeger constant $h(M)$ (see [8]).
A collection $\{ M_i \}$ of $n$-dimensional Riemannian manifolds 
(for some fixed $n$) is known as an {\sl expanding family} 
if $\inf \lambda_1(M_i) > 0$. When the Ricci
curvature of the manifolds $M_i$ is bounded from below, this definition is
equivalent to the condition that $\inf h(M_i) > 0$, by work of Cheeger [8]
and Buser [7].
The construction of explicit expanding families of
manifolds has been a major topic of research, with applications
to such diverse fields as group theory [19], lattices in Lie groups [20],
number theory [5] and coding theory [21].
More recently, it has become important to understand
the behaviour of $\lambda_1(M)$ when $M$ is a hyperbolic
3-manifold, because this has connections with the virtually
Haken conjecture, which is a major unsolved problem in
low-dimensional topology (see [15]). The Lubotzky-Sarnak conjecture [20]
proposes that any closed hyperbolic manifold has
a tower of finite covers which does not
form an expanding family. This conjecture has
major ramifications: it would imply, for example, that any arithmetic
lattice in ${\rm PSL}(2, {\Bbb C})$ has a finite index
subgroup with a non-abelian free quotient, by work of the
author, Long and Reid [16].

It is the purpose of this paper to examine $h(M)$ and $\lambda_1(M)$
for hyperbolic \break 3-manifolds $M$. Our main theorem imposes upper bounds
on $h(M)$ and $\lambda_1(M)$ in terms of data from any surgery diagram of $M$.
Applying this to the case of trivial surgery, where no solid tori are
attached, this gives upper bounds on $h(M)$ and $\lambda_1(M)$ when $M$ is the complement
of a hyperbolic link in the 3-sphere.
There will be three main applications. Firstly, we will show that a sequence
of hyperbolic alternating link complements with volumes tending to infinity
cannot form an expanding family. Secondly, we apply this result to 
establish a finiteness theorem for alternating link complements that
are congruence arithmetic 3-manifolds. Thirdly, we utilise the existence of 
expanding families of hyperbolic 3-manifolds to prove that some 3-manifolds must
have `complicated' surgery diagrams.

In order to state our main theorem, we need some terminology.
Let $D$ be a link diagram, and let $G(D)$ be the
underlying 4-valent planar graph. A {\sl bigon region} is the closure of
a complementary region of $G(D)$ with two edges in its boundary.
A {\sl twist region} is either a maximal collection of bigon
regions whose union is connected or a single crossing adjacent
to no bigon regions. The crossings of a diagram are canonically
partitioned into twist regions. The {\sl twist number} $t(D)$
is the number of twist regions. Thus, $t(D)$ is always at most
the crossing number $c(D)$ of the diagram.

\vskip 18pt
\centerline{
\epsfxsize=1.8in
\epsfbox{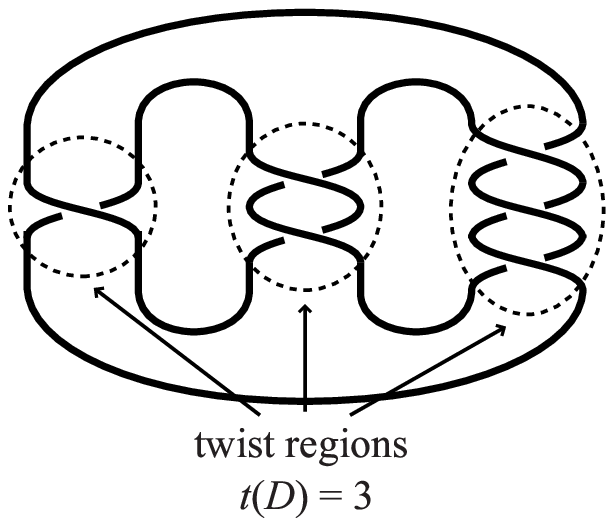}
}
\vskip 6pt
\centerline{Figure 1}

Recall that a compact orientable 3-manifold $M$ is obtained by {\sl Dehn surgery} on a link $L$
in the 3-sphere if there is a collection of properly
embedded disjoint simple closed curves in $M$, such that removing an open
regular neighbourhood of these curves results in the exterior of $L$.
Equivalently, $M$ can be constructed from the exterior of $L$ by attaching
solid tori, where the boundary of each solid torus is homeomorphically
identified with a boundary component of the exterior of $L$. The essentially 
different ways of attaching each solid torus are parametrised by an element 
of ${\Bbb Q} \cup \infty$. A {\sl rational surgery
diagram} for $M$ is a diagram for $L$, plus an assignment of an element of
${\Bbb Q} \cup \{ \infty \}$ to some components
of $L$. We permit some of the components of $L$ to be unlabelled,
in which case no solid torus is attached to the relevant component
of $\partial N(L)$.

We say that a compact 3-manifold is {\sl hyperbolic} if its interior admits
a complete finite volume hyperbolic structure. We say that a link $L$
in $S^3$ is {\sl hyperbolic} if its exterior is hyperbolic.

Our main theorem is as follows.

\noindent {\bf Theorem 1.1.} {\sl Let $D$ be a rational surgery diagram
of a compact orientable hyperbolic 3-manifold $M$. Then the Cheeger constant
$h(M)$ satisfies
$$h(M) \leq 4 \pi {(24 \sqrt 2 + 16 \sqrt 3) \sqrt {t(D)}\over {\rm Volume}(M)}
\leq 4 \pi {(24 \sqrt 2 + 16 \sqrt 3)
\sqrt {c(D)} \over {\rm Volume}(M)}.$$
}

By applying this to the case where no component of the link is filled in, we obtain
the following corollary.

\noindent {\bf Corollary 1.2.} {\sl Let $D$ be a diagram of a hyperbolic link $L$
in the 3-sphere. Then
$$h(S^3 - L) \leq 4 \pi {(24 \sqrt 2 + 16 \sqrt 3) \sqrt {t(D)} \over {\rm Volume}(S^3 - L)}
\leq 4 \pi {(24 \sqrt 2 + 16 \sqrt 3)
\sqrt {c(D)} \over {\rm Volume}(S^3 - L)}.$$
}

We now give some applications of this corollary.
If one wants to use this result to find upper bounds on $h(S^3 - L)$ from the diagram $D$, one needs
to find lower bounds on the volume of $S^3 - L$. Fortunately, such bounds are known
for two large classes of links: alternating links and highly twisted links.

We start with alternating links. If a link has an alternating diagram, then it
has one where the underlying 4-valent graph has no edge loops. This is because any such
loop may removed to produce an alternating diagram of the same link with fewer crossings.
A diagram is {\sl twist-reduced} if,
for every simple closed curve that meets the link projection transversely in
four points away from the crossings, with two points of intersection adjacent
to one crossing and the other two points of intersection adjacent to another crossing,
the simple closed curve bounds a (possibly empty) collection of bigons
arranged end to end between the crossings (see Figure 2).
If an alternating diagram is not twist-reduced, then we may produce
an alternating diagram of the same link with the same number of crossings
but smaller twist number. 

\vskip 18pt
\centerline{
\epsfxsize=3.5in
\epsfbox{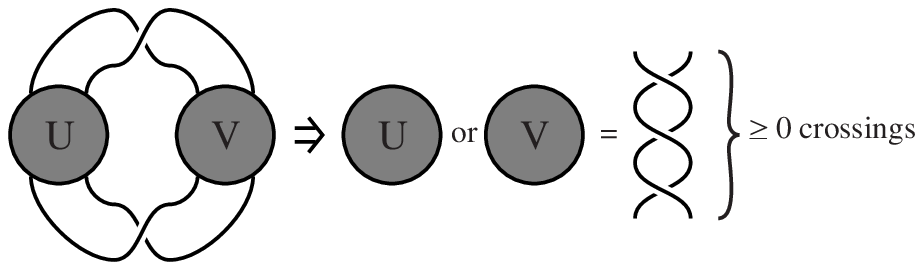}
}
\vskip 6pt
\centerline{Figure 2}

The main theorem of [13] is as follows.

\noindent {\bf Theorem 1.3.} {\sl Let $L$ be a hyperbolic link with a twist-reduced alternating
diagram $D$ having no edge loops. Then
$$v_3(t(D) - 2) \leq {\rm Volume}(S^3 - L) < 16 v_3 (t(D) - 1),$$
where $v_3 \simeq 1.0149$ is the volume of a regular hyperbolic ideal tetrahedron.}

In fact, the constants in this theorem have been improved upon. The
upper bound on the volume was reduced to $10 v_3 (t(D) - 1)$ by Agol and D. Thurston [13].
The lower bound on volume was increased to $v_8(t(D)/2 - 1)$ by Agol, Storm and W. Thurston [2],
where $v_8 \simeq 3.6639$ is the volume of a regular hyperbolic ideal octahedron. It is these improved bounds 
that we use. 

This result has been extended to another class of links
by Futer, Kalfagianni and Purcell [11]. A diagram is known as {\sl highly twisted}
if each twist region contains at least 7 crossings, and the diagram is
alternating within each twist region. We say that a link
is {\sl highly twisted} if it has a twist-reduced highly twisted diagram
having no edge loops. Their result is as follows.

\noindent {\bf Theorem 1.4.} {\sl Let $L$ be a hyperbolic link with a twist-reduced
highly twisted diagram $D$ having no edge loops. Then
$$0.70735 \ (t(D) - 1) < {\rm Volume}(S^3 - L) < 10 v_3 (t(D) - 1).$$}

By combining Corollary 1.2, Theorem 1.3 and Theorem 1.4, we obtain the following 
consequence for the geometry of these link complements.

\noindent {\bf Corollary 1.5.} {\sl Let $L$ be a hyperbolic link with an alternating
or highly twisted diagram $D$ that is twist-reduced, has no edge loops and where $t(D) > 2$. Then
$$\eqalign{
h(S^3 - L) & \leq c_1 / \sqrt {t(D)} \cr
h(S^3 - L) & \leq c_2 / \sqrt{{\rm Volume}(S^3 - L)} \cr}$$
where $c_1 \leq 1643$ and $c_2 \leq 1129$.
}

Buser's inequality [7] gives that $\lambda_1(S^3 - L) \leq 4 h(S^3 - L) + 10 (h(S^3 - L))^2$, and
hence we have the following spectral consequence.

\noindent {\bf Corollary 1.6.} {\sl Let $L$ be a hyperbolic link with an alternating
or highly twisted diagram $D$ that is twist-reduced, has no edge loops and where $t(D) > 2$. Then
$$\eqalign{
\lambda_1(S^3 - L) & \leq c_3 / \sqrt {t(D)} + c_4 / t(D) \cr
\lambda_1(S^3 - L) & \leq c_5 / \sqrt{{\rm Volume}(S^3 - L)}
+ c_6 / {\rm Volume}(S^3 - L), \cr
}$$
where $c_3 \leq 6572$, $c_4 \leq 2.7 \times 10^7$, $c_5 \leq 4516$ and $c_6 \leq 1.3 \times 10^7$.}

This implies one direction of the following.

\noindent {\bf Theorem 1.7.} {\sl A collection of alternating or highly twisted
hyperbolic link complements forms an expanding family if and only if their
volumes are bounded.}

For the other direction, suppose that an infinite collection of alternating or
highly twisted hyperbolic link complements have bounded volume, but do not form
an expanding family. Then, we may pass to a subsequence where the smallest positive eigenvalue
of the Laplacian tends to zero. We may pass to a further subsequence $M_i$ which converges in the Gromov-Hausdorff
topology to a fixed finite-volume hyperbolic 3-manifold $M_\infty$ (see Chapter E in [3] for
example). By [9], the eigenvalues $\lambda_1(M_i)$ tend to $\lambda_1(M_\infty)$, and hence are bounded away from zero, 
which is a contradiction.

This result can be interpreted in several different
ways. On the one hand, it may mean that expanding
families of hyperbolic 3-manifolds are `rare'. If so,
this would provide support for the Lubotzky-Sarnak
conjecture. On the other hand, it may mean that
alternating and highly twisted link complements are not representative
of `generic' hyperbolic 3-manifolds. 

Given the above non-expansion results for two large
classes of knot and link complements, it is natural
to speculate about all hyperbolic knot and link complements.
This leads to the following interesting question:

\noindent {\bf Question.} Does there exist a collection of hyperbolic link
complements with volumes that tend to infinity and which forms
an expanding family?

A consequence of Theorem 1.7 is the following finiteness result. We refer the reader to [1] for the definition of a
{\sl congruence} arithmetic 3-manifold. According to Theorem 5.9 of [1] (see also Corollary 1.3(a) of [6]), 
any such 3-manifold $M$ has 
$\lambda_1(M) \geq 3/4$. Also, by [4], there are only finitely many arithmetic 3-manifolds with volume
less than any given real number. Thus, we obtain the following corollary.

\noindent {\bf Corollary 1.8.} {\sl There are only finitely many alternating or highly
twisted link complements that are congruence arithmetic hyperbolic 3-manifolds.}

This leads to the following interesting question.

\noindent {\bf Question.} Are there only finitely many alternating or highly
twisted link complements that are arithmetic?

The methods of this paper cannot be immediately applied here, since there exist
arithmetic 3-manifolds with arbitrarily small Cheeger constant. However, any
arithmetic 3-manifold finitely covers a congruence arithmetic 3-manifold.
So, this may perhaps give a route to tackling this question.

We now return to surgery diagrams for 3-manifolds. It is natural
to ask how complex must a surgery diagram of a 3-manifold be.
In this direction, Constantino and D. Thurston have the following theorem [10].

\noindent {\bf Theorem 1.9.} {\sl There is a positive constant $c$ such that
any closed hyperbolic 3-manifold $M$ has a rational surgery diagram with
at most $c ({\rm Volume}(M))^2$ crossings.}

It is not hard to show that the number of crossings required in such 
a surgery diagram $D$ is at least ${\rm Volume}(M)/(10v_3)$. This is because
the exterior of the link specified by $D$ has Gromov norm at most
$10 c(D)$ by the appendix in [13]. Gromov norm does not
increase under Dehn filling [23], and so the Gromov norm of $M$ is at most
$10 c(D)$. But the Gromov norm of a hyperbolic 3-manifold $M$ is equal
to ${\rm Volume}(M)/v_3$, by [23]. This establishes the required bound.

Given this linear lower bound on crossing number, it is natural to ask whether
the quadratic upper bound of Constantino and Thurston is sharp. Here,
we show that it is, by applying Theorem 1.1 to an expanding family
of hyperbolic 3-manifolds.

\noindent {\bf Theorem 1.10.} {\sl Let $\{ M_i \}$ be an expanding
family of finite-volume orientable hyperbolic 3-manifolds. Then, there is a positive constant $c$ such
that any rational surgery diagram for $M_i$ requires at least
$c ({\rm Volume}(M_i))^2$ crossings.}

Such expanding families of 3-manifolds are known to exist. Indeed, the
following result of Long, Lubotzky and Reid [18] shows that they arise as
covering spaces of any given finite-volume hyperbolic 3-manifold.

\noindent {\bf Theorem 1.11.} {\sl Any finite-volume hyperbolic 3-manifold has an infinite
sequence of distinct finite-sheeted regular covers that forms an expanding family.}

The proof of our main result, Theorem 1.1, utilises the relationship between the following
various notions of `width':
\item{1.} The max-width of a link in the 3-sphere. This is closely related to a concept introduced
by Gabai in his proof of the Property R conjecture [12]. We give its definition below.
\item{2.} The width of an abstract graph. (By `abstract' here, we just mean that the
graph is not necessarily embedded anywhere, and if it is, then this embedding is
immaterial.) This notion was first defined by the author in [14],
where it was used to prove a group-theoretic result. We recall it in Section 2.
\item{3.} The Morse width of a planar graph. This is a new concept, which measures
the interaction between a graph embedded in the 2-sphere and Morse functions for that
sphere. We will define it in Section 3.
\item{4.} The Heegaard width of a compact orientable 3-manifold $M$. This measures the complexity of generalised
Heegaard splittings for $M$. A very closely related concept was defined by the author in [15] (and
was denoted $c_+(M)$ there). We will give its definition in Section 5.

Let us start with the definition of the max-width of a link. One considers diagrams $D$ 
for the link in ${\Bbb R}^2$. The {\sl max-width} of $D$ is the maximum, over all $t \in {\Bbb R}$, of the number of
intersections of ${\Bbb R} \times \{ t \}$ with the link projection.
The {\sl max-width} of the link is defined to be the minimal max-width of any of its
diagrams. This is a slight variant of Gabai's definition in [12].
He considered the sum of the number of intersections, over a finite collection
of representative level sets ${\Bbb R} \times \{ t \}$. We call this the
{\sl sum-width} of the diagram. (See Section 6 where the
definition of sum-width is recalled in detail.) However, here, it
is more appropriate to consider the maximum rather than the sum.

Let $D$ be a rational surgery diagram for a compact orientable hyperbolic 3-manifold $M$,
and let $G(D)$ be the underlying 4-valent planar graph. For convenience,
we assume that $D$ is a diagram in the 2-sphere rather than ${\Bbb R}^2$.
The first step in the proof of Theorem 1.1 is to bound the width of
$G(D)$ as an abstract graph. To do this, we use a famous result of Lipton and Tarjan [17] that provides
an efficient `separator' for a planar graph. This is a collection of vertices
in the graph with relatively small size, which divides the graph into
two subgraphs, each with relatively large size. Using an inductive
argument involving these separators, we can find an upper bound on
the width of $G(D)$ as an abstract graph. This is then
used to construct a Morse function on the diagram 2-sphere, so that
each level set of the function has relatively few intersections with the graph.
This extends to a Morse function on the 3-sphere, in which level
sets have relatively few intersections with the link. We arrange
that this Morse function is equivalent to the standard Morse function
on the 3-sphere, and so this provides an upper bound on the max-width of
the link, as defined above. This leads to the following result,
which is of independent interest.

\noindent {\bf Theorem 1.12.} {\sl Let $D$ be a diagram of a link $L$ in the
3-sphere. Then the max-width of $L$ is at most
$$2 + (24 \sqrt 2 + 16 \sqrt 3) \sqrt{t(D)}.$$}

Note that we are not claiming here that the max-width of $D$, or any diagram
ambient isotopic to it, is bounded by the above function.

A theorem of Thompson [22] relates the width of a knot to
its bridge number, provided that the knot is not tangle-composite.
However, some care is required here, because Thompson used the sum-width of a knot
(as defined by Gabai) rather than max-width. In Section 6, we show how these 
notions are related, and apply Theorem 1.12 to prove the following.

\vfill\eject
\noindent {\bf Theorem 1.13.} {\sl Let $D$ be a diagram of a knot $K$ in the
3-sphere that is not tangle-composite. Then the bridge number of $K$ is at most
$$1 + (12 \sqrt 2 + 8 \sqrt 3) \sqrt{t(D)}.$$}

We now return to the proof of Theorem 1.1. By Theorem 1.12, the link
specified by $D$ has a diagram with max-width at most $2 + (24 \sqrt 2 + 16 \sqrt 3) \sqrt{t(D)}$.
Using this diagram, we construct a generalised
Heegaard splitting of the link exterior, with control over the
Euler characteristic of the splitting surfaces. This extends to
a generalised Heegaard splitting of the filled-in 3-manifold $M$. Finally,
we use a result of the author from [15] which places an upper bound on the Cheeger
constant of a hyperbolic 3-manifold in terms of data from a generalised
Heegaard splitting.

\vskip 18pt
\centerline{\caps 2. The width of planar graphs}
\vskip 6pt

In this and the following two sections, we give the proof of Theorem 1.12.
The first step is as follows. Given a link diagram $D$, let $G(D)$ be
its underlying planar graph. Let $T(D)$ be the corresponding
{\sl twist graph}, which has a vertex for each twist region of $G(D)$ and an edge
for each edge of $G(D)$ not lying in a twist region. Thus, $T(D)$ is a planar
graph which is obtained by collapsing $G(D)$.
(See Figure 3.)

\vskip 18pt
\centerline{
\epsfxsize=2.8in
\epsfbox{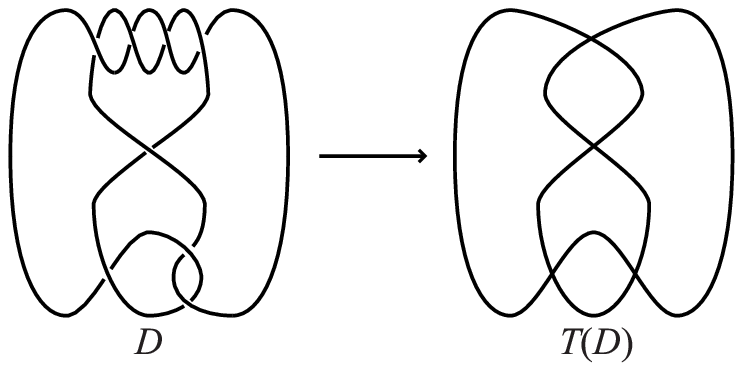}
}
\vskip 6pt
\centerline{Figure 3}

The width of a finite graph $G$ was defined in [14]. We recall
the definition here. Let $V(G)$ be the vertex set of $G$
with cardinality $v(G)$. For any subset $A$ of $V(G)$, let
$\partial A$ denote the set of edges with one endpoint in $A$
and one endpoint not in $A$. 
We consider all possible bijections $\phi \colon [1, v(G)] \cap {\Bbb N} \rightarrow V(G)$.
(Such a bijection is effectively just a total ordering on the vertices.)
The {\sl width} of $\phi$ is 
$${\rm width}(\phi) = \max_{i \in [1,v(G)] \cap {\Bbb N}} |\partial (\phi([1,i] \cap {\Bbb N}))|.$$
The {\sl width} of $G$, denoted ${\rm width}(G)$, is the minimum, over
all such $\phi$, of ${\rm width}(\phi)$.

Although it is not strictly relevant here, there is an attractive 
pictorial interpretation of the width of a finite graph $G$ (without loops and without
isolated vertices). Given a bijection
$\phi \colon [1,v(G)] \cap {\Bbb N} \rightarrow V(G)$, draw the graph in
${\Bbb R}^2$ so that the vertices have height given by $\phi^{-1}$,
and so that the height function on each edge has no critical points and so that 
the edges intersect transversely. Then,
the width of $\phi$ is equal to the maximal number of intersections,
over all $t \in {\Bbb R}$, between the graph and ${\Bbb R} \times \{ t \}$.
Hence, the width of $G$ is the minimum of this maximum, over all such
realisations of $G$ in ${\Bbb R}^2$. There are obvious analogies with
the max-width of a link defined in Section 1.

\vskip 18pt
\centerline{
\epsfxsize=2in
\epsfbox{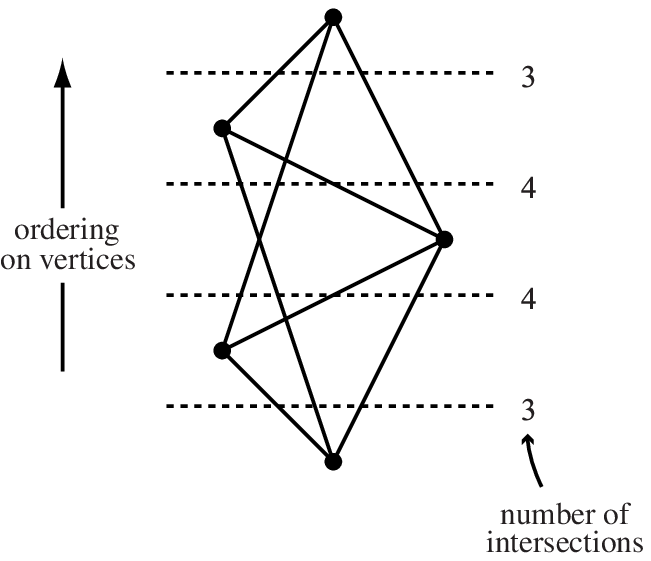}
}
\vskip 6pt
\centerline{Figure 4}

One reason for studying the width of a graph $G$ is its relation
to the graph's Cheeger constant. Recall [19] that this is defined
to be
$$h(G) = \min \left \{ {|\partial A| \over |A|}: A \subset V(G), 0 < |A| \leq v(G)/2 \right \}.$$
Setting $i$ to be $\lfloor v(G)/2 \rfloor$ in the definition of graph width,
we obtain the inequality
$$h(G) \leq {{\rm width}(G) \over \lfloor v(G)/2 \rfloor}.$$

The following result is central to this paper. It gives an upper bound
on the width of a planar graph $G$, which grows as a function of $\sqrt {v(G)}$.

\noindent {\bf Theorem 2.1.} {\sl Let $G$ be a finite planar graph,
with $v(G)$ vertices and maximal vertex degree $\Delta(G)$. Then
$${\rm width}(G) \leq {(6 \sqrt 2 + 4 \sqrt 3)} \Delta(G) \sqrt{v(G)}.$$}

The proof relies heavily on the important work of
Lipton and Tarjan on separators for planar graphs [17]. Recall that a
set of vertices $S$ in a graph $G$ {\sl separates} $G$ into
disjoint graphs $G_1$ and $G_2$ if one obtains $G_1 \cup G_2$ by
removing $S$ and all adjacent edges from $G$, and if there
is no edge joining $G_1$ to $G_2$. (We do not require either
$G_i$ to be connected.)

\noindent {\bf Theorem 2.2.} (Lipton-Tarjan [17]) {\sl Any finite planar graph
has a set of vertices $S$ which separates $G$ into $G_1$ and $G_2$,
such that
$$v(G_1) \leq (2/3) v(G), \qquad v(G_2) \leq (2/3) v(G), \qquad
|S| \leq \sqrt {8 v(G)}.$$}

\noindent {\sl Proof of Theorem 2.1.} We prove this by induction on $v(G)$.
For brevity, set $n = v(G)$. The induction starts trivially with $n = 1$.
Let us therefore prove the inductive step. Let $S$ be the set of vertices provided
by Lipton and Tarjan's theorem, which separates $G$ into $G_1$ and $G_2$.
By induction, for $j = 1$ and $2$, there is a bijection 
$\phi_j \colon [1, v(G_j)] \cap {\Bbb N} \rightarrow
V(G_j)$ with width at most
$${(6 \sqrt 2 + 4 \sqrt 3)} \Delta(G_j) \sqrt{v(G_j)} \leq 
{4 (\sqrt 3 + \sqrt 2)} \Delta(G) \sqrt n.$$
We use these to define a bijection $\phi \colon [1, n] \cap {\Bbb N} \rightarrow V(G)$, as follows.
Send the first $|S|$ integers to $S$ via some arbitrary bijection.
Then send the next $v(G_1)$ vertices to $V(G_1)$ via $\phi_1$.
(In other words, define $\phi(i) = \phi_1(i - |S|)$ in this
part of the domain.) Then send the final $v(G_2)$ vertices to $V(G_2)$ via
$\phi_2$. It is clear that 
$${\rm width}(\phi) \leq \Delta(G) |S| + \max \{ {\rm width}(\phi_1), {\rm width}(\phi_2) \}.$$
To see this, consider some integer $i$ between $1$ and $v(G)$, and 
let $A = \phi([1,i] \cap {\Bbb N} )$. If $i \leq |S|$, then any edge
in $\partial A$ is adjacent to a vertex in $S$. There are at most
$\Delta(G) |S|$ such edges. If $i$ lies
between $|S|+1$ and $|S| + v(G_1)$, any edge in $\partial A$
is either adjacent to $S$ or runs between two vertices in $G_1$.
In the latter case, this edge contributes to ${\rm width}(\phi_1)$.
A similar argument holds when $i$ lies between $|S| + v(G_1) + 1$
and $v(G)$. Thus,
$${\rm width}(G) \leq {\rm width}(\phi) \leq \Delta(G) \sqrt {8 n}  + 
{4 (\sqrt 3 + \sqrt 2)} \Delta(G) \sqrt n = {(6 \sqrt 2 + 4 \sqrt 3)} \Delta(G) \sqrt n,$$
as required. $\square$

We will apply this result to the twist graph $T(D)$ of a diagram $D$. The following can then
be used to provide an upper bound on the width of $G(D)$.

\noindent {\bf Lemma 2.3.} {\sl Let $D$ be a link diagram. Then ${\rm width}(G(D))
\leq {\rm width}(T(D))+2$.}

\noindent {\sl Proof.} Start with a total ordering $\phi$ on the vertices of $T(D)$, with minimal width.
Use this to construct a total ordering on the vertices of $G(D)$, as follows.
The first vertex of $T(D)$ corresponds to a twist region of $G(D)$. We set the vertices
of this twist region to be the first vertices in the ordering, taken in order
along the twist region. Then consider the second vertex of $T(D)$, and so on.
Let $\psi$ be the resulting bijection $[1,c(D)] \cap {\Bbb N} \rightarrow V(G(D))$. Then, for
each positive integer $i$, $\psi([1,i] \cap {\Bbb N})$ consists of the vertices in the
twist regions of $\phi([1,j] \cap {\Bbb N})$, for some positive integer $j$, plus a subset of the vertices
in a single twist region. Thus, the edges in $\partial (\psi([1,i] \cap {\Bbb N}))$ correspond
to the edges $\partial (\phi([1,j] \cap {\Bbb N}))$, plus possibly two edges in the twist
region. Hence, the width of $\psi$ is at most ${\rm width}(\phi) + 2$,
which equals ${\rm width}(T(D)) + 2$. $\square$

\vskip 18pt
\centerline{\caps 3. The Morse width of planar graphs}
\vskip 6pt

Let $G$ be a finite graph embedded in the 2-sphere, in which each
vertex has valence at most $4$.
In this section, our goal is to construct a Morse function $f \colon S^2 \rightarrow {\Bbb R}$,
starting from a total ordering of the vertices of $G$.

Let $f \colon S^2 \rightarrow {\Bbb R}$
be a Morse function. We say that $f$ is {\sl generic} with respect to $G$ if its critical points
have distinct values, and the vertices of $G$ have distinct values. However,
we allow the possibility that a critical point and a vertex take the same
value; indeed, we allow vertices to be critical points of $f$.
The {\sl width} of $G$ with respect to $f$ is defined to be $\max \{ |f^{-1}(t) \cap G|
: t \in {\Bbb R} \}$. The {\sl Morse width}
of the embedded graph $G$ is the minimal width of $G$ with respect to $f$, over all
generic Morse functions $f$. 

\noindent {\bf Proposition 3.1.} {\sl Let $G$ be a finite
graph embedded in $S^2$, in which each
vertex has valence at most $4$. Then the Morse width of $G$ is equal to its
width as an abstract graph.}

\noindent {\sl Proof.} We first show that 
the width of $G$ is at most its Morse width.
Consider a generic Morse function $f$ on the 2-sphere. The vertices are ordered according
to their values under $f$, and so this defines a bijection $\phi \colon [1,v(G)] \cap {\Bbb N}\rightarrow V(G)$.
Let $t_i$ be the value under $f$ of the $i^{\rm th}$ vertex.
If $t \in (t_i, t_{i+1})$, then clearly, any edge 
in $\partial (\phi([1,i]\cap {\Bbb N}))$ must go through $f^{-1}(t)$. Hence,
$|f^{-1}(t) \cap G|$ is at least $|\partial (\phi([1,i]\cap {\Bbb N}))|$.
The width of $\phi$ is therefore at most
the width of $G$ with respect to $f$, which establishes the required
inequality.

To prove the inequality in the other direction,
suppose that we are given a bijection $\phi \colon [1, v(G)] \cap {\Bbb N}\rightarrow
V(G)$. We will extend $\phi^{-1}$ to a Morse function $f$ 
on $S^2$ with the same width. First define
$f$ to be monotonic on the edges. Since each vertex has valence
at most $4$, we may extend $f$ to a regular neighbourhood of $G$.
(See Figure 5 for some examples.) Then extend $f$ 
over each complementary region of $G$ to create a generic Morse function.
Clearly, the width of $G$ with respect to $f$ is equal 
to the width of $\phi$, as required. $\square$

\vskip 6pt
\centerline{
\epsfxsize=3.4in
\epsfbox{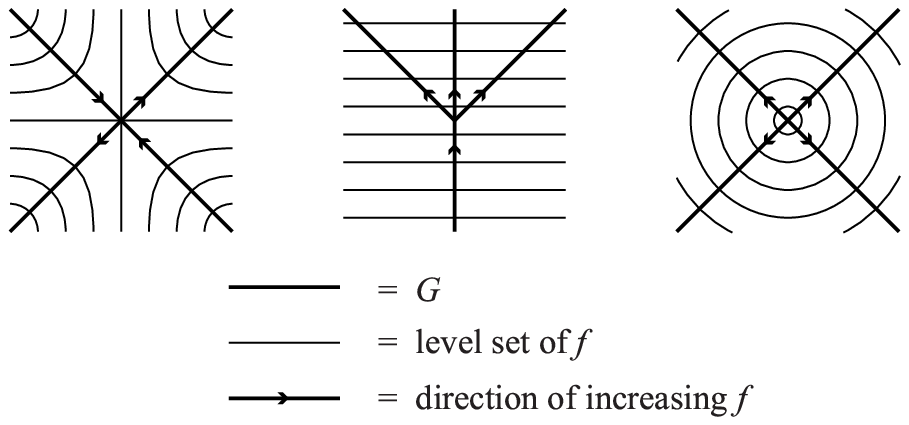}
}
\vskip 6pt
\centerline{Figure 5}

An example of an embedded planar graph $G$ is given in Figure 6.
There, a Morse function on $S^2$ is constructed from a bijection
$[1,v(G)] \cap {\Bbb N} \rightarrow V(G)$, using the recipe in
the above proof.

\vskip 18pt
\centerline{
\epsfxsize=3.2in
\epsfbox{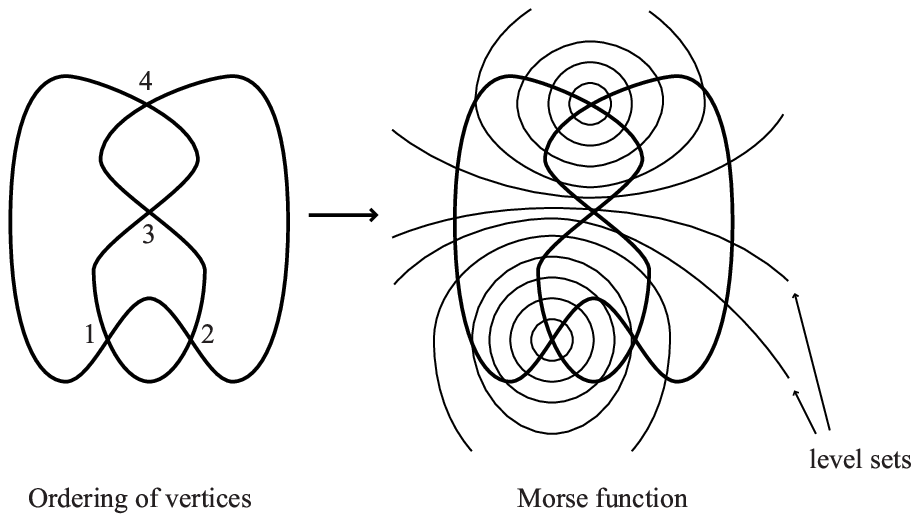}
}
\vskip 6pt
\centerline{Figure 6}

\vskip 18pt
\centerline{\caps 4. The max-width of links}
\vskip 6pt

In Section 1, we gave the definition of the max-width of a link $L$ in the 3-sphere.
We now give an interpretation of max-width which is slightly more topological.
Let $f \colon S^3 \rightarrow {\Bbb R}$ be the standard Morse function with
a single maximum, a single minimum and no other critical points. Let $L$ be
a link embedded in $S^3$ so that $f|L$ is a Morse function. The {\sl max-width}
of this embedding is $\max \{ |f^{-1}(t) \cap L| : t \in {\Bbb R} \}$. It is clear that the
minimum max-width over all links ambient isotopic to $L$ is equal to the
max-width of $L$, as defined in Section 1.

We now give a bound on the max-width of $L$ in terms of Morse width, as defined in Section 3.

\noindent {\bf Proposition 4.1.} {\sl Let $D$ be a diagram of a link $L$. Then
the max-width of $L$ is at most the Morse width of the embedded graph $G(D)$.}

\noindent {\sl Proof.} Let $f \colon S^2 \rightarrow {\Bbb R}$ be a Morse
function that is generic with respect to $G(D)$ and which realises the Morse width of $G(D)$. From this, we
will construct a Morse function $F \colon S^3 \rightarrow {\Bbb R}$,
equivalent to the standard Morse function, such that the max-width of $L$ with respect
to $F$ is equal to the width of $G(D)$ with respect to $f$.
Some care is required here, because $f$ may have many critical points, whereas
we require $F$ to have just two.

Let us first rescale $f$ so that its image lies in $(-1,1)$.

Let $S^2 \times [-2,2]$ be a regular neighbourhood of the equatorial
sphere in $S^3$. We may pick an embedding of $L$ in $S^2 \times [-1,1]$
as specified by the diagram $D$. More specifically, we may arrange
that the projection of $L$ onto the first factor $S^2$ of the product
is equal to the link projection in $D$, and that the behaviour of $L$
near the inverse image of the double points coincides with the crossing information
of $D$.

We will now specify a foliation of $S^2 \times [-2,2]$ by 2-spheres. The leaves will be
indexed by $t \in [-1,1]$, and will be denoted by $S^2_t$. This
foliation will be equivalent to the product foliation, and so
the leaves $S^2_t$ can be viewed as the level sets of a smooth function
$F \colon S^2 \times [-2,2] \rightarrow [-1,1]$ without critical points. We can then
extend this to a Morse function on $S^3$ that is equivalent to the
standard Morse function.

For values of $t$ away from small neighbourhoods of the critical values, 
we define $S^2_t$ initially to be
$$
( f^{-1}[t,1] \times \{ t-1 \} )
\cup
( f^{-1}(t) \times [t-1, t+1] )
\cup 
( f^{-1}[-1,t] \times \{ t+1 \} ).
$$
See Figure 7 for an example of such an $S^2_t$. It is clear that for regular values
of $t$, $S^2_t$ is a 2-sphere. For it is obtained by cutting $S^2$ along $f^{-1}(t)$
to give $f^{-1}[-1,t]$ and $f^{-1}[t,1]$ and then inserting annuli between the corresponding
boundary components of these surfaces. It is also straightforward to check
that these spheres are pairwise disjoint.

\vskip 18pt
\centerline{
\epsfxsize=3.5in
\epsfbox{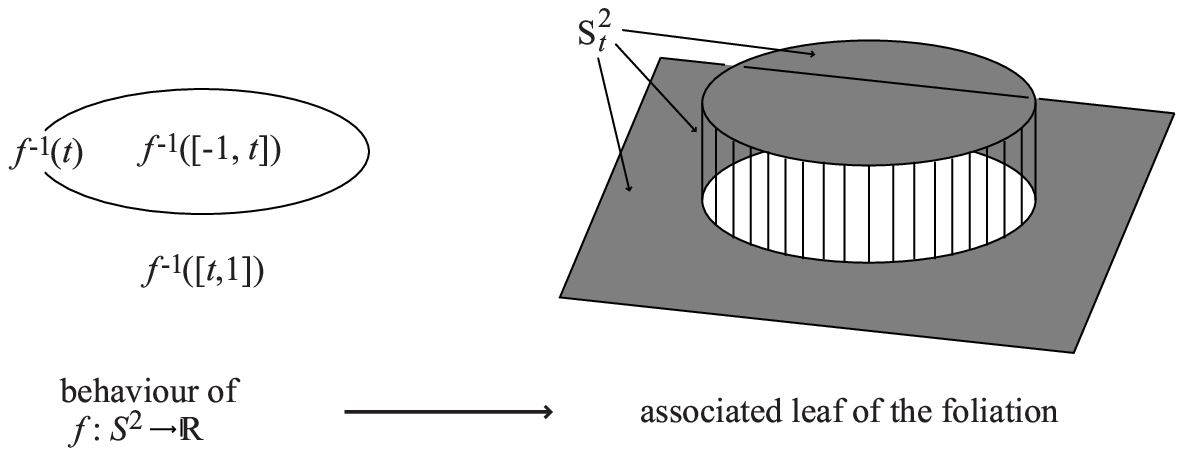}
}
\vskip 6pt
\centerline{Figure 7}

Near the critical values of $f$, we need to modify this definition slightly,
in order that $S^2_t$ is genuinely a sphere. (For example, if a local maximum of $f$
has value $t$, then $S^2_t$ as defined above would be the union of a sphere
and an interval.) However, this modification can clearly
be done. We also need to modify the entire foliation so that it is smooth.
Again, this is clearly possible.

An example is given in Figure 8. In the first diagram, the level sets of the
Morse function are shown. In the second diagram, a cross-section of the preliminary `foliation'
is depicted. In the third diagram, the modification that makes the foliation
smooth is shown.

\vskip 18pt
\centerline{
\epsfxsize=4.5in
\epsfbox{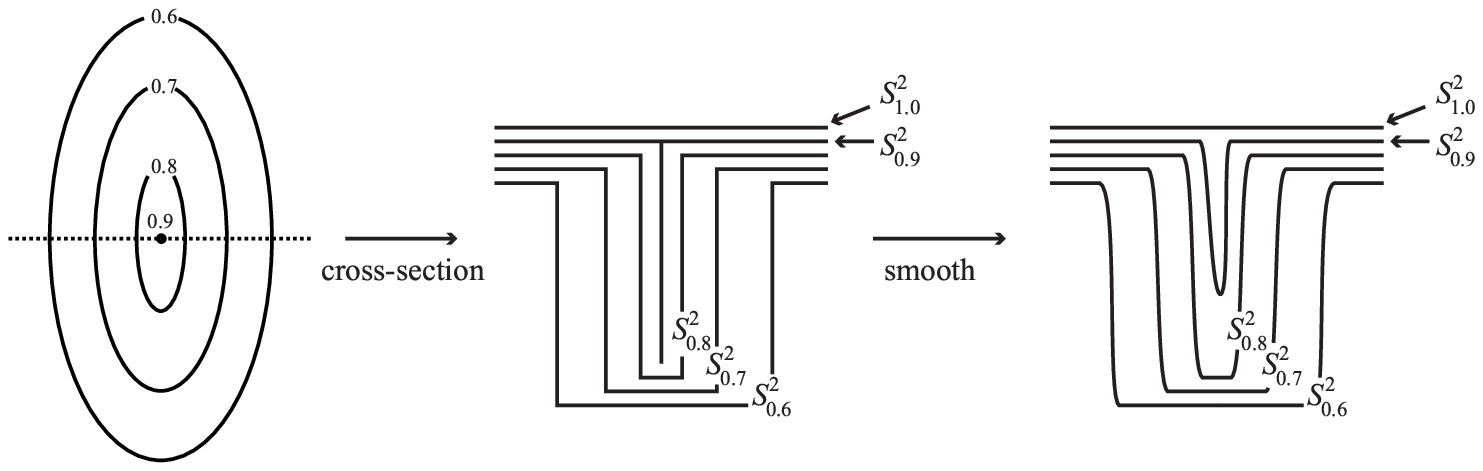}
}
\vskip 6pt
\centerline{Figure 8}

For any value of $t$ in $[-1,1]$, the number of
intersections between $S^2_t$ and $L$
is equal to the number of intersections
between $f^{-1}(t)$ and $G(D)$. Hence, the max-width of $L$ 
with respect to $F$ is indeed equal to the width of $G(D)$ with
respect to $f$, as required. $\square$

We can now complete the proof of Theorem 1.12.

\noindent {\sl Proof.} Let $D$ be a diagram of a link $L$ in the 3-sphere.
By Proposition 4.1, the max-width of $L$ is at most the Morse width of $G(D)$.
By Proposition 3.1, this is equal to the width of $G(D)$ as an abstract graph.
By Lemma 2.3, this is at most ${\rm width}(T(D)) + 2$.
This is at most  ${(24 \sqrt 2 + 16 \sqrt 3)} \sqrt{t(D)} + 2$,
by Theorem 2.1. $\square$

\vskip 18pt
\centerline{\caps 5. The Heegaard width of 3-manifolds}
\vskip 6pt

Let $M$ be a compact orientable 3-manifold. In this section, we will 
study an invariant of $M$ known as its {\sl Heegaard width}.
This is defined in terms of the complexity of surfaces in 
generalised Heegaard splittings of $M$.

Recall that, for a compact connected surface $S$, $\chi_-(S)$
is $\max \{ 0, -\chi(S) \}$. When $S$ is a compact,
possibly disconnected surface, $\chi_-(S)$ is defined to be
the sum of $\chi_-(S_i)$, as $S_i$ ranges over all the
components of $S$.

Let $\{ C_1, \dots, C_n \}$ be a generalised Heegaard splitting for $M$. Thus, 
each $C_i$ is a disjoint union of compression bodies; $M$ is union of these compression
bodies; their interiors are disjoint; their boundary components are either
properly embedded or components of $\partial M$; and, for each relevant odd integer $i$,
$\partial_+ C_i = \partial_+ C_{i+1}$ and $\partial_- C_i  \cap {\rm int}(M) = \partial_- C_{i-1}  \cap {\rm int}(M)$.
We define the {\sl width} of the splitting to be
$\max_i \chi_-(\partial_+ C_i)$. The {\sl Heegaard width} of
$M$ is the minimal width of any generalised Heegaard splitting.
We denote it by ${\rm Heeg}\hbox{-}{\rm width}(M)$.

In this section, we will prove the following.

\noindent {\bf Proposition 5.1.} {\sl Let $L$ be a link in the 3-sphere.
Then the Heegaard width of the exterior of $L$ is at most ${\rm max\hbox{-}width}(L) - 2$.}

\noindent {\sl Proof.} Let $f \colon S^3 \rightarrow {\Bbb R}$ be the
standard Morse function, and place $L$ so that it has minimal max-width
with respect to $f$. We may assume that $f|L$ is a Morse function and that the critical points of $f|L$
occur at distinct heights. Let $S_1, \dots, S_n$ be a minimal collection of
level sets of $f$, with increasing values under $f$, so that each $S_i$
avoids the critical points of $f|L$ and so that 
between $S_i$ and $S_{i+1}$ only local minima or only local maxima of $f|L$ appear.
In other words, if we view the critical points of $f|L$ as a sequence
of local minima, followed by a sequence of local maxima, followed by a sequence 
of local minima, and so on, then the surfaces $S_i$ occur between the
 local minima and the local maxima, and between the local maxima and the local minima.
From the surfaces $S_i$,
we construct closed surfaces $F_i$ in the exterior of $L$, as follows.

Around $L$, place $n$ parallel tori $T_1, \dots, T_n$, where $T_1$ is closest to $\partial N(L)$,
$T_2$ is adjacent to $T_1$, and so on. For each $i$, the intersection $T_i \cap S_i$
is a collection of simple closed curves, which bound discs in $S_i$.
Remove the interiors of these discs, and attach the parts of $T_i$ that
lie above $S_i$. Let $F_i$ be the resulting surface. (See Figure 9.) 

We claim that 
the surfaces $F_1, \dots, F_n$ 
separate the exterior of $L$ into compression bodies $C_1, \dots, C_{n+1}$. To see this,
consider the submanifold $C_{i+1}$ between $F_i$ and $F_{i+1}$. This is a
modification of the region between $S_i$ and $S_{i+1}$. Suppose first 
that this region only contains local maxima of $L$ and
that $i \not= n$. For each local maximum of $L$ between $S_i$ and $S_{i+1}$,
there is a disc $D$ such that $L \cap D$ is an arc in $\partial D$
running from $S_i$ up to the local maximum of $L$ and back down to $S_i$, and such that the
remainder of $\partial D$ is an arc in $S_i$. We may arrange that
these discs are all disjoint and that the intersection between each disc 
and $C_{i+1}$ is a compression disc for
$F_i$. If we cut $C_{i+1}$ along these discs, the resulting manifold is
homeomorphic to $F_{i+1} \times I$ (see Figure 10).

\vskip 18pt
\centerline{
\epsfxsize=3.5in
\epsfbox{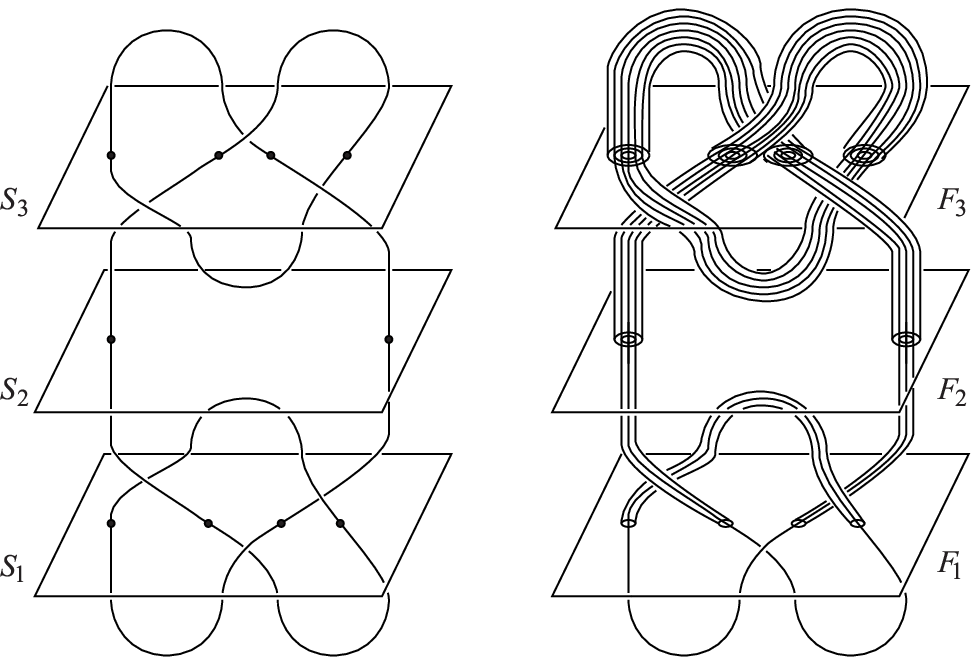}
}
\vskip 6pt
\centerline{Figure 9}

\vskip 18pt
\centerline{
\epsfxsize=3in
\epsfbox{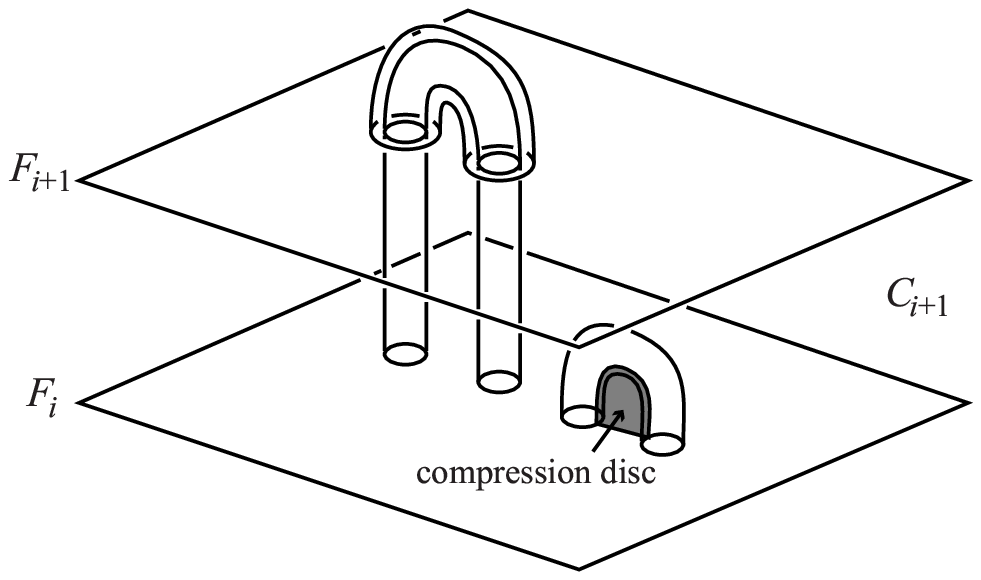}
}
\vskip 6pt
\centerline{Figure 10}

Suppose now that the region between $S_i$ and $S_{i+1}$ contains only local
minima of $L$ and that $i \not= 1$. As previously, for each local minimum, there 
is a disc $D$ such that $L \cap D$ is an arc in $\partial D$
running from $S_{i+1}$ down to the local minimum of $L$ and back up to $S_{i+1}$, and such that the
remainder of $\partial D$ is an arc in $S_{i+1}$. Let $N$ be a big regular neighbourhood of
$D$, and let
$D'$ be $\partial N \cap C_{i+1}$. Then, provided $N$ was big enough,
$D'$ is a disc properly embedded in $C_{i+1}$ with boundary
in $F_{i+1}$ (see Figure 11). We may arrange that
these discs $D'$ are all disjoint. If we cut $C_{i+1}$ along
these discs, the result is a collection of copies of $T^2 \times I$
and a copy of $F_i \times I$. Thus, $C_{i+1}$ is a compression body.
Similarly, $C_1$ and $C_{n+1}$ are compression bodies.

These compression bodies form a generalised Heegaard splitting for
the exterior of $L$. Its width is 
$\max_i \{ |S_i \cap L| - 2\}$, which equals ${\rm max\hbox{-}width}(L) - 2$.
$\square$

\vskip 18pt
\centerline{
\epsfxsize=3.7in
\epsfbox{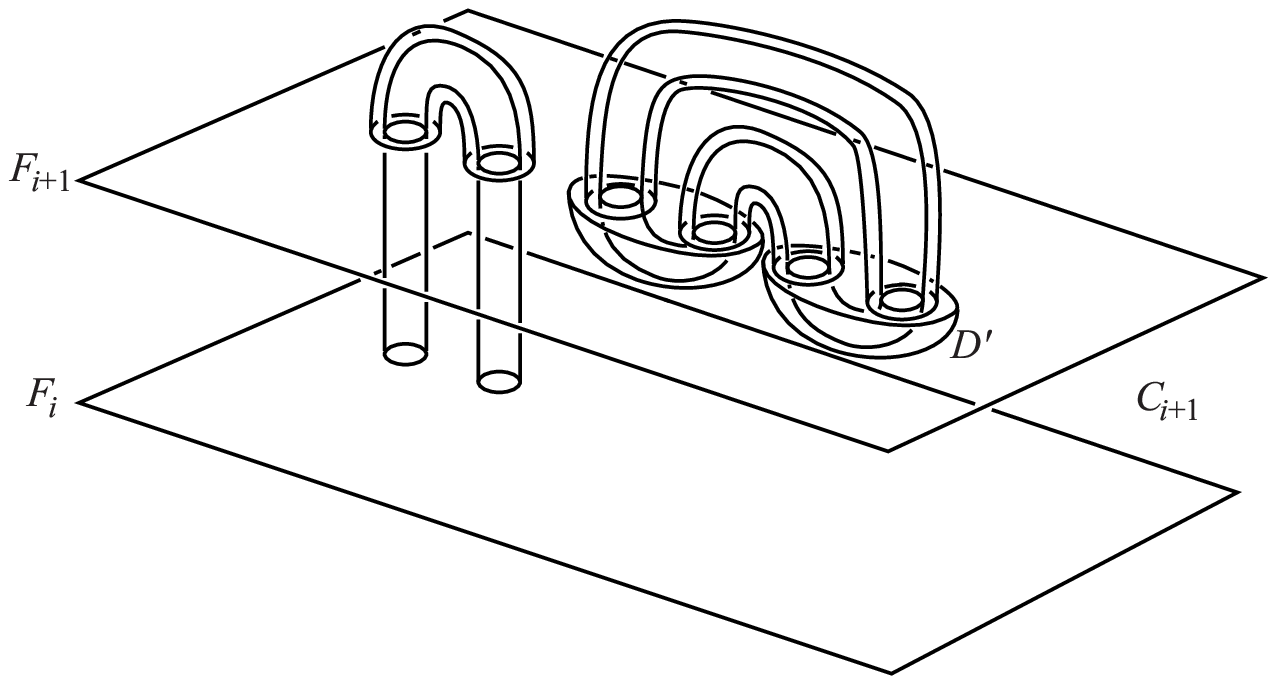}
}
\vskip 6pt
\centerline{Figure 11}

The key relationship between the Heegaard width of a finite-volume hyperbolic manifold and
its Cheeger constant is contained in the following result. The following
is a slight variant of Theorem 4.1 in [15], and it has essentially the same
proof, which we will not repeat here.

\noindent {\bf Theorem 5.2.} {\sl Let $M$ be a compact orientable finite-volume hyperbolic 3-manifold.
Then
$$h(M) \leq 4 \pi {{\rm Heeg}\hbox{-}{\rm width}(M)
\over {\rm Volume}(M)}.$$}

\noindent {\sl Proof of Theorem 1.1.}
Let $D$ be a rational surgery diagram
of a compact orientable hyperbolic 3-manifold $M$.
Let $L$ be the link defined by $D$.
By Proposition 5.1, the Heegaard width of the exterior of $L$
is at most ${\rm max\hbox{-}width}(L) - 2$. By Theorem 1.12,
this at most ${(24 \sqrt 2 + 16 \sqrt 3)} \sqrt{t(D)}$.
Now, the Heegaard width of a compact orientable 3-manifold does not increase under Dehn filling, 
since a generalised Heegaard splitting for the unfilled manifold
becomes a generalised Heegaard splitting for the filled-in one. Hence,
$${\rm Heeg}\hbox{-}{\rm width}(M)
\leq {(24 \sqrt 2 + 16 \sqrt 3)} \sqrt{t(D)}.$$
Applying Theorem 5.2 gives the required upper bound on $h(M)$.
$\square$

\vskip18pt
\centerline{\caps 6. Other notions of link width}
\vskip 6pt

In this section, our goal is to prove Theorem 1.13, which places
an upper bound on the bridge number of any tangle-prime knot, in
terms of the twist number of any of its diagrams. In order to do this,
we need to investigate yet more notions of the width of a link.
We have already defined the max-width of a link in Section 1.
We now give some variants of this.

\noindent {\caps Sum-width}

This is the original version of width
which Gabai introduced in his proof of the
Property R conjecture. Let $D$ be a diagram for a link $L$ in ${\Bbb R}^2$.
Projection onto the second co-ordinate of ${\Bbb R}^2$ gives
a height function. Let us suppose that the restriction to the link
of this height function is Morse, and the critical points of the link have
distinct values. These values divide ${\Bbb R}$ into open intervals,
and for each value of $t$ in one of these intervals, the number
of intersection between the link projection and ${\Bbb R} \times \{ t \}$
is constant. Let us define the {\sl sum-width} of $D$ to be
the sum of the number of intersection points, over each of these
intervals. The {\sl sum-width} of $L$ is the minimal sum-width
of any diagram for $L$.

\noindent {\caps Lex-width}

Although sum-width is very useful and elegant,
one loses a certain amount of information when performing
the summation. Instead, one can combine the intersection numbers of the intervals 
into a multi-set (that is, a set where repetitions
are retained), called the {\sl lex-width} of $D$. One compares two multi-sets using
the usual lexicographical ordering. Here,
one re-orders each multi-set into a descending
sequence. One compares the first two integers
of each multi-set. If they are the same, one
passes to the second two, and so on. The
minimal lex-width (with respect to this ordering)
over all diagrams $D$ for $L$ is known as the
{\sl lex-width} of $L$. Unlike sum-width and max-width,
this is not a single integer, but is a multi-set
of integers. The following is trivial.

\noindent {\bf Lemma 6.1.} {\sl The largest integer in the multi-set
${\rm lex\hbox{-}width}(L)$ is ${\rm max\hbox{-}width}(L)$.}

\noindent {\sl Proof.} Let $D$ be a diagram for $L$ with minimal lex-width. 
The max-width of $D$ is the maximum integer in
the multi-set ${\rm lex\hbox{-}width}(D)$, which equals ${\rm lex\hbox{-}width}(L)$.
Thus, the max-width of $L$ is at most this maximum.

Conversely, let $D'$ be a diagram for $L$ with minimal max-width.
The lex-width of $D'$ is an upper bound for the lex-width of $L$.
But when comparing multi-sets, one compares their largest integers
first. So, the maximal integer in ${\rm lex\hbox{-}width}(L)$ is
at most the max-width of $D'$, which equals the max-width of $L$. $\square$

A diagram for a link $L$ is a {\sl bridge diagram} if projection onto the
second factor of ${\Bbb R}^2$ restricts to a Morse function on $L$, in which all the local
maxima occur above all the local minima. The {\sl bridge number} $b(L)$ is
the minimal number of local maxima in any bridge diagram.
Hence, we trivially have that ${\rm max\hbox{-}width}(L) \leq
2 b(L)$. We will now see that this is,
in fact, an equality for tangle-prime knots.

It is a theorem of Thompson [22] that, when $K$ is
a tangle-prime knot, then every diagram for $K$ of minimal sum-width is
a bridge diagram. However, exactly
the same argument gives the following
related result.

\noindent {\bf Theorem 6.2.} {\sl If a
knot $K$ is tangle-prime, then every diagram for $K$ of minimal lex-width is
a bridge diagram.}

\noindent {\bf Corollary 6.3.} {\sl
If $K$ is a tangle-prime knot,
then ${\rm max\hbox{-}width}(K) = 2 b(K)$.}

\noindent {\sl Proof.} The largest integer
in the multi-set ${\rm lex\hbox{-}width}(K)$ is
${\rm max\hbox{-}width}(K)$. But, in a bridge
diagram, the maximal number of intersections between ${\Bbb R} \times \{ t \}$
and the link projection is at least $2 b(K)$. 
Thus, by Theorem 6.2, ${\rm max\hbox{-}width}(K) \geq 2 b(K)$.
Since the opposite inequality always
holds, this must be an equality. $\square$

We can now prove Theorem 1.13.

\noindent {\sl Proof.} Let $D$ be a diagram of a knot $K$ in the
3-sphere that is not tangle-composite. By Corollary 6.3, the bridge number of $K$ 
is equal to half the max-width of $K$. By Theorem 1.12, this is at most
$1 + (12 \sqrt 2 + 8 \sqrt 3) \sqrt{t(D)}$. $\square$

\vskip 18pt
\centerline{\caps References}
\vskip 6pt

\item{1.} {\caps I. Agol, M. Belolipetsky, P. Storm, K. Whyte}, {\sl
Finiteness of arithmetic hyperbolic reflection groups},
arXiv:math.GT/0612132

\item{2.} {\caps I. Agol, P. Storm, W. Thurston}, {\sl
Lower bounds on volumes of hyperbolic Haken 3-manifolds},
J. Amer. Math. Soc. 20 (2007) 1053--1077.

\item{3.} {\caps R. Benedetti, C. Petronio,} 
{\sl Lectures on hyperbolic geometry.} Universitext. Springer-Verlag, Berlin, 1992

\item{4.} {\caps A. Borel}, {\sl Commensurability classes and volumes of hyperbolic
3-manifolds}, Ann. Scuola Norm. Sup. Pisa 8 (1981) 1--33.

\item{5.} {\caps J. Bourgain, A. Gamburd, P. Sarnak,}
{\sl Sieving and expanders.} C. R. Math. Acad. Sci. Paris 343 (2006), no. 3, 155--159.

\item{6.} {\caps M. Burger, P. Sarnak}, {\sl Ramanujan duals II}, Invent.
Math. 106 (1991) 1-11.

\item{7.} {\caps P. Buser,} {\sl A note on the isoperimetric constant}, Ann. Sci.
Ecole Norm. Sup. 15 (1982) 213--230.

\item{8.} {\caps J. Cheeger}, {\sl A lower bound for the smallest eigenvalue of
the Laplacian}, Problems in analysis (Papers dedicated to Salomon Bochner, 1969),
pp 195--199, Princeton University Press.

\item{9.} {\caps B. Colbois, G. Courtois}, {\sl Convergence de vari\'et\'es et convergence
du spectre du Laplacien}, Ann. Sci. Ecole Norm. Sup. 24 (1991) 507--518.

\item{10.} {\caps F. Constantino, D. Thurston}, {\sl 
$3$-manifolds efficiently bound $4$-manifolds}, J. Topology 1 (2008) 703--745.

\item{11.} {\caps D. Futer, E. Kalfagianni, J. Purcell,} {\sl 
Dehn filling, volume and the Jones polynomial}, J. Differential Geom.
78 (2008) 429--464.

\item{12.} {\caps D. Gabai,} {\sl Foliations
and the topology of 3-manifolds, III},
J. Differential Geom. 26 (1987) 479--536.

\item{13.} {\caps M. Lackenby,} {\sl The volume of hyperbolic alternating
link complements}, Proc. London Math. Soc. 88 (2004) 204--224 (with
an appendix by I. Agol and D. Thurston)

\item{14.} {\caps M. Lackenby}, {\sl A characterisation of large
finitely presented groups}, J. Algebra 287 (2005) 458--473.

\item{15.} {\caps M. Lackenby}, {\sl Heegaard splittings, the virtually Haken
conjecture and Property $(\tau)$}, Invent. Math. 164 (2006) 317--359.

\item{16.} {\caps M. Lackenby, D. Long, A. Reid}, {\sl Covering spaces of arithmetic
3-orbifolds}, Int. Math. Res. Not. (2008)

\item{17.} {\caps R. Lipton, R. Tarjan}, {\sl A separator theorem for
planar graphs}, SIAM J. Applied Mathematics 36 (1979) 177--189.

\item{18.} {\caps D. D. Long, A. Lubotzky, A. W. Reid}, {\sl
Heegaard genus and Property $(\tau)$ for hyperbolic 3-manifolds,}
J. Topology 1 (2008) 152--158.

\item{19.} {\caps A. Lubotzky}, {\sl Discrete Groups,
Expanding Graphs and Invariant Measures}, \hfill\break Progress
in Math. 125 (1994)

\item{20.} {\caps A. Lubotzky}, {\sl Eigenvalues of the 
Laplacian, the first Betti number and the congruence subgroup
problem,} Ann. Math. 144 (1996) 441--452.

\item{21.} {\caps M. Sipser, D. Spielman}, {\sl Expander codes},
IEEE Trans. on Information Theory, 42 (1996) 1710--1722.

\item{22.} {\caps A. Thompson,} {\sl Thin position and bridge number for knots 
in the $3$-sphere.} Topology 36 (1997) 505--507.

\item{23.} {\caps W. Thurston,} {\sl The geometry and topology of 3-manifolds,}
Lecture notes (1980).

\vskip 12pt
\+ Mathematical Institute, University of Oxford, \cr
\+ 24-29 St Giles', Oxford OX1 3LB, United Kingdom. \cr

\end